\documentclass[12pt,reqno]{amsart}

\usepackage[all]{xy}
\usepackage{amsmath}
\usepackage{amsfonts}
\usepackage{amssymb}
\usepackage[all]{xy}
\usepackage{mathrsfs}
\usepackage{bbding}
\usepackage{txfonts}
\usepackage{amscd}

\usepackage[shortlabels]{enumitem}
\usepackage{ifpdf}
\ifpdf
\usepackage[colorlinks,final,backref=page,hyperindex]{hyperref}
\else
\usepackage[colorlinks,final,backref=page,hyperindex,hypertex]{hyperref}
\fi
\usepackage{tikz}
\usepackage[active]{srcltx}

\def\W{\mathcal{SD}}

\def\M{M}
\def\Z{\mathbb{Z}}
\def\N{\mathbb{N}}

\def\C{\mathbb{C}}

\numberwithin{equation}{section}
\newtheorem{theo}{Theorem}[section]
\newtheorem{defi}[theo]{Definition}
\newtheorem{coro}[theo]{Corollary}
\newtheorem{lemm}[theo]{Lemma}
\newtheorem{prop}[theo]{Proposition}

\newtheorem{case}{Case}

\newtheorem{rema}[theo]{Remark}

\newtheorem{remark}[theo]{Remark}

\newcommand{\nc}{\newcommand}

\nc{\tred}[1]{\textcolor{red}{#1}}
\nc{\tblue}[1]{\textcolor{blue}{#1}}
\nc{\tgreen}[1]{\textcolor{green}{#1}}
\nc{\tpurple}[1]{\textcolor{purple}{#1}}
\nc{\btred}[1]{\textcolor{red}{\bf #1}}
\nc{\btblue}[1]{\textcolor{blue}{\bf #1}}
\nc{\btgreen}[1]{\textcolor{green}{\bf #1}}
\nc{\btpurple}[1]{\textcolor{purple}{\bf #1}}

\nc{\yf}[1]{\textcolor{red}{Yufeng:#1}}
\nc{\ld}[1]{\textcolor{blue}{Liu Dong: #1}}
\nc{\hb}[1]{\textcolor{purple}{Haibo: #1}}

\begin{document}

\title{Irreducible Modules for Super-Virasoro Algebras from Algebraic D-Modules }

\author{Haibo Chen}
\address{Department of Mathematics, Jimei University, Xiamen, Fujian 361021, China}
\email{hypo1025@163.com}

\author{Xiansheng Dai}
\address{School of Mathematical Sciences, Guizhou Normal University,
Guiyang 550001,  China}
\email{daisheng158@126.com}

\author{Dong Liu}
\address{Department of Mathematics, Huzhou University, Zhejiang Huzhou, 313000, China}
\email{liudong@zjhu.edu.cn}

\author[Yufeng Pei]{Yufeng Pei}
\address{Department of Mathematics, Shanghai Normal University, Guilin Road 100,
	Shanghai 200234, China}\email{pei@shnu.edu.cn}

\vspace{-5mm}

\begin{abstract}

In this paper, we  introduce a new family of functors from the category of modules for the Weyl algebra to the category of  modules for the super-Virasoro algebras. The properties of these functors are investigated, with an emphasis on  irreducibility preservation and natural isomorphisms. By utilizing these functors, we  recover some old irreducible super-Virasoro modules, including those from the irreducible intermediate series as well as irreducible $U(\mathfrak{h})$-free modules. Additionally, we provide several families of new irreducible super-Virasoro modules via our constructed functors.

\end{abstract}



\maketitle
\vspace{-10mm}

\section{Introduction}

The super-Virasoro algebras are crucial infinite-dimensional Lie superalgebras in theoretical physics. They serve as supersymmetric extensions of the ordinary Virasoro algebra and consist of two sectors: the Neveu-Schwarz sector \cite{NS} and the Ramond sector \cite{R}.
The study of string theory and conformal field theory relies heavily on representations of the super-Virasoro algebra, as they describe the quantum states of a supersymmetric theory. The highest weight modules for the super-Virasoro algebras have been extensively studied in literature (see \cite{FQS,KW, MR}, etc.).  In \cite{IK}, researchers investigated Verma modules' structure for the super-Virasoro algebras, constructing Bernstein-Gelfand-Gelfand type resolutions for irreducible highest weight modules. Additionally, Fock modules over super-Virasoro algebras were analyzed in \cite{IK1}.

Recently, there have been several advancements in the study of representation theory. For example, the classification of simple Harish-Chandra modules, i.e., irreducible modules with
finite-dimensional weight spaces, for the Neveu-Schwarz and Ramond algebras has been achieved in \cite{CLL,CL,S}.
Every such module is either an irreducible highest or an irreducible lowest weight module, or an irreducible module of the intermediate series. Furthermore, in \cite{LPX}, the concept of Whittaker modules for these algebra was introduced and analyzed. In addition to this, \cite{C,LPX1} explored the study of simple restricted modules that generalize both highest weight and Whittaker modules for the Neveu-Schwarz and Ramond algebras.  In \cite{YYX1}, the authors presented a family of non-weight modules for the super-Virasoro algebra, called $U(\mathfrak{h})$-free modules. They demonstrated that these modules provide a complete classification of free $U(\mathfrak{h})$-modules of rank 1 over the Ramond algebra. Moreover, they established that the same family of modules also constitutes a complete classification of free $U(\mathfrak{h})$-modules of rank 2 over the Neveu-Schwarz algebra.

The theory of Weyl algebras, which are also known as algebras of differential operators, is an interesting and important subject.
The study of algebraic D-modules, which are modules over Weyl algebras \cite{CC}, has a significant impact on various fields.
Used the twisting technique,  L\"{u} and Zhao in \cite{LZ}  created lots of new irreducible Virasoro modules by twisting irreducible modules over the Weyl algebra. For additional insights and further references, please refer to \cite{LMG} and \cite{PSTZ}.

This paper aims to extend the twisting technique introduced in \cite{LZ} to the super-Virasoro algebra.  We present a method for constructing new irreducible modules for the super-Virasoro algebras, based on irreducible modules for the Weyl algebra,
We start by defining a functor $S$ that maps modules over the Weyl algebra $\mathcal{D}$ to modules over the Weyl superalgebra $\mathcal{SD}$. We prove that this functor preserves irreducibility and establishes a one-to-one correspondence between irreducible modules over the two algebras. Next, we introduce a family of functors $F_{\epsilon,b}$ that operate on modules for the Weyl superalgebra and yield modules for the centerless super-Virasoro algebra $\mathfrak{g}[\epsilon]$. These functors depend on parameters $\epsilon$ (which can be either 0 or $\frac12$) and $b\in \C$. By composing these functors with $S$, we obtain a new family of functors $H_{\epsilon, b}$ that operate on modules for the Weyl algebra and produce modules for the super-Virasoro algebra.  Using Block's classification of irreducible Weyl modules in \cite{B}, together with our newly defined functors, we establish that these functors preserve irreducibility for all values of $b$ except for a few cases which are also clearly described.  We also determine necessary and sufficient conditions under which two such functors are naturally isomorphic.
Our main results are presented below.

\vspace{0.2cm}

\noindent {\bf Main theorem  } \ [Theorem \ref{main2}  and Theorem \ref{th4333}] \ {\it
Suppose that $\epsilon\in\{0,\frac{1}{2}\},b,b_1,b_2\in\mathbb{C}$. Then
there exists a functor
$H_{\epsilon,b}:{\rm Mod}\,\mathcal{D}\to {\rm Mod}\,\mathfrak{g}[\epsilon]$ with
$
M\mapsto H_{\epsilon,b}(M),
$
where $H_{\epsilon,b}(M)=M\oplus\overline{M}$ and is defined by the following actions:
\begin{eqnarray*}
L_m\cdot v&=&-t^{m}\left(D+(m-2\epsilon)b+\frac{m+2\epsilon}{2}\theta\partial_\theta\right)v, \\
G_p\cdot v &=&t^{p-\epsilon}\left(\theta D+2(p-\epsilon)b\theta-t^{2\epsilon}\partial_\theta\right)v,
\end{eqnarray*}
 for $m\in\mathbb{Z}$, $p\in\mathbb{Z}+\epsilon$,  and $v\in M\oplus\overline{M}$. Furthermore, we show
that
\begin{itemize}
    \item[{\rm(i)}]  $H_{\epsilon,b}$ preserves irreducible modules if and only if $b\notin\{0,\frac{1}{2}\}$.
    \item[{\rm (ii)}] $H_{\epsilon,b_1}\cong H_{\epsilon,b_2}$  if and only if  $b_1=b_2$.
\end{itemize}
}

 Applying these results, we  recover some old irreducible super-Virasoro modules, including those from the irreducible intermediate series \cite{S} as well as irreducible $U(\mathfrak{h})$-free modules \cite{YYX1}. Finally, we provide several examples of new irreducible super-Virasoro modules via our constructed functors.

Our approach for constructing modules for the super-Virasoro algebra can also be applied to the $N=2$ superconformal algebras. For more information on these modules, please refer to our forthcoming paper \cite{CLP}.

This paper is organized as follows: Section 2 provides a brief review of the definitions of super-Virasoro algebras, Weyl algebra, and Weyl superalgebra. In Section 3, we introduce a family of functors that map modules for the Weyl algebra to modules for the super-Virasoro algebras. Section 4 investigates properties of irreducibility-preserving and natural isomorphisms for these functors. Finally, in Section 5, we utilize these functors to construct several new irreducible modules for the super-Virasoro algebras.

Throughout this paper, $\C,\C^*,\Z, \Z^*, \N$ and $\Z_+$ will denote the sets of complex numbers, nonzero complex numbers, integers, nonzero integers, nonnegative integers, and positive integers, respectively.  All vector spaces, vector algebras, and modules will be considered to be over the field of complex numbers $\C$, and all modules for Lie superalgebras will be $\Z_2$-graded.

\section{Preliminaries}

In this section, we recall some notations, and collect the known facts about Lie superalgebras, the super-Virasoro  algebras, and the Weyl (super)algebra.

\subsection{} Suppose that  $L$ is a Lie (super)algebra (or an associative (super)algebra), we shall denote the category  of $L$-modules by  ${\rm Mod}\; L$.
  We will also denote the parity change functor by $\Pi$, which is defined as follows: $\Pi(V)= V_{\bar 1}\oplus V_{\bar 0}$, where $(\Pi(V))_{\bar0}= V_{\bar 1}$ and $(\Pi(V))_{\bar1}=V_{\bar0}$, where $V=V_{\bar 0}\oplus V_{\bar1}$ is an $L$-module.  For a Lie superalgebra ${L}$, we will use ${U}(L)$ to denote its universal enveloping algebra, and use $\delta_{i,j}$ to denote the Kronecker delta.
\begin{defi}
Let $M$ be a module of an  associative (super)algebra   $L$.
If   $0\neq v\in M$ and some $x\in L$, there exists $n\in\Z_+$ such that $x^nv=0$,  then we call that  the action
of   $x$  on $v$  is nilpotent.
\end{defi}

 \begin{defi}
Suppose that  $L$ is an associative (super)algebra (or a Lie (super)algebra) and $E$ is a subspace of $L$. Let $M$ be an  $L$-module. If there exists   $0\neq x\in E$ such that $xv=0$ for some   $0\neq v\in M$,
 the  $M$  is called {\it $E$-torsion}; otherwise $M$ is called {\it $E$-torsion-free}.
\end{defi}

\begin{defi}
Let $L$ be a Lie (super)algebra with a Cartan subalgebra  $H$. An $L$-module $M$ is called a weight module if the action of $H$ on $M$ is
diagonalizable.
\end{defi}

\subsection{}

The following notion of the super-Virasoro algebra is a natural   supersymmetric generalization of the Virasoro algebra, which was first discovered by Neveu-Schwarz and Ramond, respectively.

\begin{defi}[\cite{NS},\cite{R}]For $\epsilon\in\{0,\frac{1}{2}\}$,
 the super-Virasoro algebra $\hat{\mathfrak{g}}[\epsilon]$ is an infinite dimensional Lie superalgebra whose even part is spanned by $\left\{L_n,C\mid n \in \mathbb{Z}\right\}$ and odd part is spanned by $\left\{G_p \mid p \in \mathbb{Z}+\epsilon\right\}$ subject to the following relations
\begin{eqnarray*}\label{def1.1}
&&[L_m,L_n] = (m-n)L_{m+n}+\frac{m^3-m}{12}\delta_{m+n, 0}C,\\
&&[G_p,G_q]= 2L_{p+q}+\frac13(p^2-\frac14)\delta_{p+q, 0}C,\\
&&[L_m,G_p]= \left(\frac{m}{2}-p\right)G_{m+p},\quad [\hat{\mathfrak{g}}[\epsilon],C]=0,
\end{eqnarray*}
where $m,n\in\mathbb{Z}$ and $p,q\in\mathbb{Z}+\epsilon$. $\hat{\mathfrak{g}}[\frac{1}{2}]$ and $\hat{\mathfrak{g}}[0]$
are called the Neveu-Schwarz algebra and the Ramond  algebra,
respectively.
\end{defi}

The even component of $\hat{\mathfrak{g}}[\epsilon]$ is isomorphic to the Virasoro algebra, which is denoted by ${\rm Vir}$.  A $\hat{\mathfrak g}[\epsilon]$-module $W$ is said to have a central charge of $c$ if the operator ${C}$ acts on $W$ as a complex scalar $c$. In this paper, we exclusively focus on $\hat{\mathfrak{g}}[\epsilon]$-modules with central charge zero, for $\epsilon=0, \frac12$. Let
$$
\mathfrak{g}[\epsilon]=\hat{\mathfrak{g}}[\epsilon]/\mathfrak{z},
$$
where $\mathfrak{z}=\C C$.
The Cartan subalgebra of $\mathfrak{g}[\epsilon]$ can be represented by $\mathfrak{h}=\mathbb{C}L_0$.

\subsection{}
Denote by
$
\mathcal{A}=\mathbb{C}\left[t, t^{-1}, \theta\right]
$
the associative superalgebra of Laurent polynomials associated to an even, $t$, and an odd, $\theta$, formal variable. Assume that $\theta^2=0$ and $\theta t=t \theta$.   Let  $ \partial_t=\frac{d}{dt}$, $D=t\frac{d}{dt}$ and $\partial_\theta =\frac{d}{d\theta}$.

\begin{defi}
The Weyl algebra $\mathcal{D}$ is the associative algebra of regular differential operators on the circle $S^1$generated by  $t^{\pm1}$ and  $\partial_t$,  subject to the relations  $\partial_t t - t \partial_t = 1.$
\end{defi}

 All irreducible $\mathcal{D}$-modules  have been classified by R.Block as follows:

\begin{theo}[\cite{B,LZ}]\label{211}
	Let $M$ be an  irreducible  $\mathcal{D}$-module.
	\begin{itemize}
		\item[\rm(1)]  If $M$ is $\C[t^{\pm1}]$-torsion-free, then $M\cong \mathcal{D}/(\mathcal{D}\cap\big(\C(t)[ D]\tau)\big)$ for some irreducible element   $\tau$    in the associative algebra $\C(t)[ D]$;
			\item[\rm(2)]  If $M$ is $\C[t^{\pm1}]$-torsion, then $M\cong \Omega(\lambda)$, where $\Omega(\lambda)$ is an irreducible $\mathcal{D}$-module with
 $\Omega(\lambda)=\C [ D]$, for any $\lambda\in\C^*$, and
  $$
  t^m D^n=\lambda^m( D-m)^n,\quad  D D^n= D^{n+1},\quad\forall n\in\N,m\in\Z.
  $$
	\end{itemize}
\end{theo}

\begin{defi}[\cite{CW}]
The Weyl superalgebra $\mathcal{SD}$ is  the associative superalgebra of regular differential operators
on the supercircle $S^{1 \mid 1}$  generated by  $t^{\pm1}$ and $\theta$, $\partial_t$, $\partial_{\theta}$, subject to the relations
$$\partial_t t - t \partial_t = 1, \quad \partial_{\theta} \theta + \theta \partial_{\theta} = 1, \quad \partial_t\theta = \partial_{\theta} t = 0.$$
\end{defi}

Note that the following elements
$$
t^k D^l, t^k D^l \theta \partial_\theta,t^k D^l \theta, t^k D^l \partial_\theta,\quad  k \in \mathbb{Z}, l \in \mathbb{N}
$$
form a linear basis of $\mathcal{S D}$, where   the odd elements $\theta$ and $\partial_\theta$ generate the four-dimensional Clifford superalgebra. Clearly $\mathcal{SD} $ is isomorphic to the tensor algebra of the Weyl algebra $\mathcal{D}$ and this Clifford superalgebra.



\section{Constructions of functors}
In this section, we shall construct a family of functors from the category of modules for the Weyl algebra $\mathcal{D}$ to the category of modules for the super-Virasoro algebra  $\mathfrak{g}[\epsilon]$.

\subsection{From $\mathcal{D}$-modules to $\mathcal{SD}$-modules}

\begin{lemm}\label{lemm3.1}
Let $M$ be a $\mathcal{D}$-module, and let $\overline M$ denote a copy of $M$ which we can turn into a $\mathcal{D}$-module using the operation
\begin{equation}
x\cdot\overline{v}=\overline{xv},\quad \forall x\in \mathcal{D},\overline v\in \overline M.\label{S1}
\end{equation}
We define
\begin{equation}
\partial_\theta \cdot v=0,\quad \theta\cdot v=\overline v,\quad \partial_\theta\cdot\overline{v}=v,\quad \theta\cdot\overline{v}=0,\quad \forall v\in M. \label{S2}
\end{equation}
Then $S(M)=M\oplus \overline M$ is an $\mathcal{SD}$-module
\end{lemm}
\begin{proof}
 It is easy to verify that $S(M)=M\oplus \overline M$ is an $\mathcal{SD}$-module.
\end{proof}

\begin{prop}\label{S}Let  ${\rm Mod}\,\mathcal{D}$ denote the category of $\mathcal{D}$-modules
 and ${\rm Mod}\,\mathcal{SD}$ denote the category of $\mathcal{SD}$-modules.
There exists a functor $S:{\rm Mod}\,\mathcal{D}\to {\rm Mod}\,\mathcal{SD}$ with $M\mapsto S(M)$, where  $S(M)=M\oplus \overline M$ is defined by  \eqref{S1} and \eqref{S2}.  Moreover, the functor $S$ preserves irreducible objects.
\end{prop}
\begin{proof}
 Suppose $M$ is an irreducible $\mathcal{D}$-module. Using \eqref{S1} and \eqref{S2}, we can define a functor $S$ such that $S(M)=M\oplus \overline M$. It follows that $S(M)$ is an irreducible $\mathcal{SD}$-module.
\end{proof}

\begin{prop}\label{I}
There is a one-to-one correspondence between the irreducible $\mathcal{SD}$-modules and irreducible $\mathcal{D}$-modules.
\end{prop}

\begin{proof}
Proposition \ref{S} states that if $M$ is an irreducible $\mathcal{D}$-module, then $S(M)$ is an irreducible $\mathcal{SD}$-module. Conversely, if $V$ is an irreducible $\mathcal{SD}$-module, then $V$ is also a $\mathcal{D}$-module. Let $M$ be an irreducible $\mathcal{D}$-submodule of $V$. Using $\mathcal{SD}=\C[\theta,\partial_\theta]\mathcal D$, we have $V=\C[\theta,\partial_\theta]M$. Note that $\partial_\theta M$ is also a $\mathcal{D}$-module. If $\partial_\theta M\neq 0$, we can replace $M$ with $\partial_\theta M$ and assume that $\partial_\theta M=0$. Therefore, $V=M\oplus \overline{M}=S(M)$.
\end{proof}


\subsection{From $\mathcal{SD}$-modules to $\mathfrak{g}[\epsilon]$-modules }\label{sections3.3}

\begin{prop}\label{lek22}
Let  $V$ be an $\mathcal{SD}$-module and $b\in\mathbb{C}$. Then $V$ is  a $\mathfrak{g}[\epsilon]$-module  with
\begin{equation}\label{LI2.1} L_m\cdot v=-t^{m}\left(D+(m-2\epsilon)b+\frac{m+2\epsilon}{2}\theta\partial_\theta\right)v, \end{equation}
 \begin{equation}\label{LI2.2} G_p\cdot v =t^{p-\epsilon}\left(\theta D+2(p-\epsilon)b\theta-t^{2\epsilon}\partial_\theta\right)v,
 \end{equation}
 for $m\in\mathbb{Z}$, $p\in\mathbb{Z}+\epsilon$,  and $v\in V$. We denote  this $\mathfrak{g}[\epsilon]$-module  by $F_{\epsilon,b}(V)$.
\end{prop}

\begin{proof}
Note that ${\mathfrak{g}}[\epsilon]$ can be embedded into the Lie superalgebra $\mathcal{SD}^-$, which is associated with the associative superalgebra $\mathcal{SD}$,
by the following  linear map $\Phi: {\mathfrak{g}}[\epsilon]\hookrightarrow \mathcal{SD}^-$ with
\begin{eqnarray*}
	&&L_m\longmapsto -t^{m}D-\frac{m+2\epsilon}{2}t^m\theta\partial_{\theta},\\
&&G_p\longmapsto t^{p-\epsilon}\theta D-t^{p+\epsilon}\partial_\theta,
\end{eqnarray*}
where $m\in\Z,p\in\mathbb{Z}+\epsilon.$  Then $\mathcal{A}$ is a natural $\mathfrak{g}[\epsilon]$-module.  Consider the extended  super-virasoro algebra $\tilde{\mathfrak{g}}[{\epsilon}]= \mathfrak{g}[\epsilon]\ltimes \mathcal A$   with the
brackets
$$[x,f]=-(-1)^{|x||f|}[f,x]=x(f),\ [f,g]=0,\quad \forall x\in \mathfrak{g}[\epsilon],f,g\in \mathcal A.$$
For $\epsilon\in\{0,\frac{1}{2}\},b\in \mathbb{C}$,  we define a linear map $\sigma_{\epsilon,b}: \tilde{\mathfrak{g}}[{\epsilon}]\to \tilde{\mathfrak{g}}[{\epsilon}]$ by
\begin{eqnarray*}\label{LI4455e2.1}
\sigma_{\epsilon,b}(L_m)&=&-t^{m}D-\frac{m+2\epsilon}{2}t^m\theta\partial_{\theta}-(m-2\epsilon)bt^{m},\\
\label{LI4455e2.2} \sigma_{\epsilon,b}(G_p)&=&t^{p-\epsilon}\theta D-t^{p+\epsilon}\partial_\theta+2b(p-\epsilon)t^{p-\epsilon}\theta,\\
\sigma_{\epsilon,b}(t^m)&=& {t^m},\\
\sigma_{\epsilon,b}(t^m\theta)&=& t^m\theta,
\end{eqnarray*}
for $m\in\Z,p\in\mathbb{Z}+\epsilon$.  It is  clear that $\sigma_{\epsilon,b}$ is an isomorphism of $\tilde{\mathfrak{g}}[{\epsilon}]$.

 Let $V$ be an $\mathcal{SD}$-module. Then $V$ is a natural $\mathcal{SD}^-$-module. It is clear that $V$
can be seen as a $\tilde{\mathfrak{g}}[{\epsilon}]$-module via the isomophism $\Phi$. We can define a new action of $\tilde{\mathfrak{g}}[{\epsilon}]$ on $V$ as follows
\begin{eqnarray*}\label{3w3e}
x\cdot v = \sigma_{\epsilon,b}(x)v,\quad\forall x\in\tilde{\mathfrak{g}}[{\epsilon}],v\in V.
\end{eqnarray*}
We denote this new module by
$F_{\epsilon,b}(V)$, is called the  twisted module of $V$ by $\sigma_{\epsilon,b}$. It is clear that the module $F_{\epsilon,b}(V)$
can be seen as a $\mathfrak{g}[\epsilon]$-module by restriction to the  subalgebra $\mathfrak{g}[\epsilon]$ of $\tilde{\mathfrak{g}}[{\epsilon}]$.
It follows that
\begin{eqnarray*}\label{wq44222.1}
&& L_m \cdot v=\sigma_{\epsilon,b}(L_m )v=-t^{m}\left(D+(m-2\epsilon)b+\frac{m+2\epsilon}{2}\theta\partial_\theta\right)v,
\\
&&G_p\cdot v=\sigma_{\epsilon,b}(G_p)v=t^{p-\epsilon}\left(\theta D+2(p-\epsilon)b\theta-t^{2\epsilon}\partial_\theta\right)v
\end{eqnarray*}
for $ m\in\Z,p\in\mathbb{Z}+\epsilon,v\in V$.
\end{proof}

\begin{theo}\label{Fun}
For $b\in\mathbb{C}$, there exists a functor
$$F_{\epsilon,b}:{\rm Mod}\,\mathcal{SD}\to {\rm Mod}\,\mathfrak{g}[\epsilon]$$ with
$
V\mapsto F_{\epsilon,b}(V),
$
where $F_{\epsilon,b}(V)$ is defined by (\ref{LI2.1})-(\ref{LI2.2}).
\end{theo}
\begin{proof}
 The sequence of Lie superalgebras
 $$
\mathfrak{g}[\epsilon] \hookrightarrow \tilde{\mathfrak{g}}[{\epsilon}] \stackrel{\sigma_{\epsilon,b}}{\longrightarrow}  \tilde{\mathfrak{g}}[{\epsilon}] \hookrightarrow \mathcal{SD}^-,
 $$
induces a sequence of  categories of modules for the Lie superalgebras
$$
{\rm Mod}\,\mathcal{SD}^-\to {\rm Mod}\,\tilde{\mathfrak{g}}[{\epsilon}]\to {\rm Mod}\,\tilde{\mathfrak{g}}[{\epsilon}]\to {\rm Mod}\,\mathfrak{g}[\epsilon].
$$
It is clear that each $\mathcal{SD}$ serves as a $\mathcal{SD}^-$-module.
It follows from Proposition \ref{lek22} that the functor $F_{\epsilon,b}$ can be written as the composition of the following four functors as follows:
$$
F_{\epsilon,b}:{\rm Mod}\,\mathcal{SD}\to {\rm Mod}\,\mathcal{SD}^-\to {\rm Mod}\,\tilde{\mathfrak{g}}[{\epsilon}]\to {\rm Mod}\,\tilde{\mathfrak{g}}[{\epsilon}]\to {\rm Mod}\,\mathfrak{g}[\epsilon].
$$
\end{proof}

The main result in this section can be derived by utilizing Proposition \ref{S} and Theorem \ref{Fun}.

\begin{theo}\label{M}
For $b\in\mathbb{C}$, there exists a functor
$$H_{\epsilon,b}:{\rm Mod}\,\mathcal{D}\to {\rm Mod}\,\mathfrak{g}[\epsilon]$$ with
$
M\mapsto H_{\epsilon,b}(M),
$
where $H_{\epsilon,b}(M)=(F_{\epsilon,b}\circ S)(M)=F_{\epsilon,b}(S(M))$.
\end{theo}

\section{Properties of functors $H_{\epsilon,b}$
}\label{ir4}

This section aims to show that, with the exception of $b\in\{0,\frac{1}{2}\}$, the functor $H_{\epsilon,b}$ (as stated in Theorem \ref{M}) preserves irreducible objects. Additionally, we will investigate the natural isomorphisms between the functors $H_{\epsilon,b}$ for $b\in \C$.

\subsection{Preservation of irreducibility}

\begin{lemm}\label{lemm2.3} Let $V$ be an $\mathcal{SD}$-module,  and $b\in\mathbb C, b\neq0,\frac{1}{2}$. For any $p\in\mathbb{Z}+\epsilon$  and $d\in\Z^*$, we define $T_{p,d}$  as follows:
\begin{eqnarray*}
T_{p,d} &=& \frac{1}{2d^2}(L_{-d}\cdot G_{p+d}+L_{d}\cdot G_{p-d}-2L_0\cdot G_{p})
\end{eqnarray*}
  in $U(\mathfrak{g}[\epsilon])$. Then, for all   $v \in V$, we have:
$$\frac{1}{b(1-2b)}T_{p,d}\cdot v= -t^{p-\epsilon}\theta v.$$

\end{lemm}
 \begin{proof}

For  $n\in\Z,p\in\mathbb{Z}+\epsilon$ and $v\in V$, we have
\begin{eqnarray*}
&&L_n\cdot G_{p-n}\cdot v
\\&=&-\big(t^{n} D+(n-2\epsilon)b t^{n}+\frac{n+2\epsilon}{2}t^{n}\theta\partial_\theta\big)\big(t^{p-n-\epsilon}\theta D+2(p-n-\epsilon)b t^{p-n-\epsilon}\theta-t^{p-n+\epsilon}\partial_\theta\big)v
\\&=&-\big(t^{p-\epsilon}\theta D^2+(p+2b(p-2\epsilon))t^{p-\epsilon}\theta D
-t^{p+\epsilon} D\partial_\theta
\\&&+(p-\epsilon)(p-2b\epsilon)2b t^{p-\epsilon}\theta
-(p+(1-2b)\epsilon)t^{p+\epsilon}\partial_\theta
\\&&+\big((-\frac{1}{2}-b)t^{p-\epsilon}\theta D+b(2b (p+\epsilon)+(\epsilon-3p))t^{p-\epsilon}\theta+(1-b)t^{p+\epsilon}\partial_\theta\big)n
\\&&+b(1-2b)n^2t^{p-\epsilon}\theta\big) v.
\end{eqnarray*}
 Choosing $n=-d,0,d$ in above equation, respectively, we obtain that

 $$\frac{1}{ b(1-2b)}T_{p,d}\cdot v= -t^{p-\epsilon}\theta v.$$
\end{proof}

\begin{lemm}\label{th36544433}
Let  $V$ be an irreducible module over the Weyl superalgebra $\W$ and $b\in \mathbb{C}$, then the $\mathfrak{g}[\epsilon]$-module $F_{\epsilon,b}(V)$ is irreducible if  $b\notin\{0,\frac{1}{2}\}$.
\end{lemm}
\begin{proof}

Suppose that  $b\neq0,\frac{1}{2}.$ Choose a nonzero element $v$ from  $V_{\bar0}$. For any   $p\in\mathbb{Z}+\epsilon$, and nonzero integer $d$, we can use the linear operator  $T_{p,d}$  defined in Lemma \ref{lemm2.3} to obtain the following results:
\begin{eqnarray}\label{322}
&&-\frac{1}{b(1-2b)}T_{p,d}\cdot v=\overline{t^{p-\epsilon}v}\in U(\mathfrak{g}[\epsilon])v,
\\&&\frac{1}{b(1-2b)}G_{-\epsilon}\cdot T_{p,d}\cdot v=t^{p-\epsilon} v\in U(\mathfrak{g}[\epsilon])v.
\end{eqnarray}
For $n\in\N$, one has
\begin{eqnarray}\label{266}
(-L_0+2b\epsilon)^n\cdot v= D^nv\in U(\mathfrak{g}[\epsilon])v.
\end{eqnarray}
From \eqref{322}-\eqref{266}, we can conclude that $U(\mathfrak{g}[\epsilon])v=F_{\epsilon,b}(V).$ This implies that $F_{\epsilon,b}(V)$ is an irreducible $\mathfrak{g}[\epsilon]$-module.
\end{proof}

\begin{lemm}\label{lemm200}
Let $M$ be an irreducible module over the Weyl algebra $\mathcal{D}$.
 \begin{itemize}
\item[\rm (i)]  The $\mathfrak{g}[\epsilon]$-module $H_{\epsilon,0}(M)$ is reducible  if and only if   $M$ is isomorphic to the $\mathcal{D}$-module $\C[t^{\pm1}];$
\item[{\rm (ii)}]
If $M=\C[t^{\pm1}]$, then $H_{\epsilon,0}(\C[t^{\pm1}])$ has a maximum submodule $\C$. Consequently, the factor module $H_{\epsilon,0}(\C[t^{\pm1}])/\C$ is an irreducible $\mathfrak{g}[\epsilon]$-module.
 \end{itemize}
\end{lemm}
\begin{proof}
{\rm (i)} Let $\overline{v}$ be any nonzero element in  $\overline{M}$. Now we consider the action of $D$ on $\overline{v}$.
\begin{case}\label{case1}
The
action of $D$ on $\overline{v}$ is nilpotent.
\end{case}
It is clear that there exists a nonzero element $\overline{v}\in\overline{M}$ such that $\overline{Dv}=0$. This implies that $\overline{D({t}^k{v})}=k\overline{t^k{v}}$. Therefore, we can conclude that $V$ is isomorphic to the $\W$-module $\C[t^{\pm1}]\oplus\overline{\C[{t}^{\pm1}]}$.  It can be easily verified that the $\mathfrak{g}[\epsilon]$-module $H_{\epsilon,0}(M)$ has a one-dimensional submodule $\C$, and hence, it is reducible.

\begin{case}
The
action of $ D$ on $\overline{v}$ is non-nilpotent.
\end{case}
We can choose a nonzero $
\overline{w}=\overline{D^kv}\in \overline{M},
$ where $k\in\Z_+$.  For $m\in\Z,d\in\Z^*$, we define
$$
Q_{m,d}=\frac{2}{d^2}(L_{-d}\cdot L_{m+d}+L_{d}\cdot L_{m-d}-2
L_0\cdot L_{m}).$$
By the similar calculation in Lemma \ref{lemm2.3}, we can get   $$Q_{m,d}\cdot\overline{w}=\overline{t^m w}\in U(\mathfrak{g}[\epsilon]){v},
$$ where $m\in\Z, d\in\Z^*$.  Using this, one gets
$$-G_{-\epsilon}\cdot(\overline{t^mw})=t^mw\in U(\mathfrak{g}[\epsilon]){v}.$$
 For $n\in\N$, one has
$$
(-L_0-\epsilon)^n\cdot \overline{w}= \overline{D^nw}\in U(\mathfrak{g}[\epsilon]){v}.
$$
  Thus, $H_{\epsilon,0}(M)=U(\mathfrak{g}[\epsilon]){v}$. In other words,  $H_{\epsilon,0}(M)$ is an irreducible $\mathfrak{g}[\epsilon]$-module  in this case.

{\rm (ii)}
It is clear that $\big(\C[t^{\pm1}]\oplus\overline{\C[{t}^{\pm1}]}\big)/ \C=\mathrm{span}_{\C}\{t^n, \overline{t^k} \mid n \in\Z^*,k\in\Z\}$. For
$n\in\Z^*,k\in\Z$, $m\in\Z$, we have
\begin{eqnarray*}
&&L_m\cdot t^n=-(t^m D)t^n=-nt^{m+n},\\
&&L_m\cdot\overline{t^k}=-\left(t^m D+\frac{m+2\epsilon}{2}t^m\theta\partial_\theta\right)\overline{t^k}=-\left(k+\frac{m+2\epsilon}{2}\right)\overline{t^{m+k}},
\\
&&G_p\cdot t^n=(t^{p-\epsilon}\theta D-t^{p+\epsilon}\partial_\theta)t^n=n\overline{t^{p-\epsilon+n}},\\
&&G_p\cdot \overline{t^k}=(t^{p-\epsilon}\theta D-t^{p+\epsilon}\partial_\theta)\overline{t^k}=-t^{p+\epsilon+k}.
\end{eqnarray*}
This shows  that $\big(\C[t^{\pm1}]\oplus\overline{\C[{t}^{\pm1}]}\big)/ \C$ is an irreducible $\mathfrak{g}[\epsilon]$-module.
\end{proof}

Suppose that $M$ is an irreducible $\mathcal{D}$-module, we define
$$
{H_{\epsilon,0}'(M)}=\begin{cases}
&H_{\epsilon,0}(M)/\C,\quad \text{if}\quad  M\cong \C[t^{\pm1}];\\
&H_{\epsilon,0}(M),\quad\quad\text{otherwise}.
\end{cases}
$$

\begin{lemm}\label{lemm201}
Suppose that $M$ is an irreducible $\mathcal{D}$-module.
\begin{itemize}
    \item[\rm (i)]  The module $H_{\epsilon,\frac{1}{2}}(M)$ is an irreducible $\mathfrak{g}[\epsilon]$-module if and only if $D(M) = M$.
    \item[\rm (ii)] If $D(M) = M$, then $H_{\epsilon,\frac{1}{2}}(M)\cong\Pi\circ {H'_{\epsilon,0}(M)}$, where $\Pi$ is the  parity change functor for $\mathfrak{g}[\epsilon]$-modules.

\end{itemize}

\end{lemm}
\begin{proof}
{\rm(i)}
Consider   $M\oplus\overline{M}$ as a $\Z_2$-graded modules over $\W$. Let $v\in M$ and $\overline{v} \in \overline{M}$. For
any $m\in\Z$ and $p\in\Z+\epsilon$,  $\mathfrak{g}[\epsilon]$-module $H_{\epsilon,\frac{1}{2}}(M)$ is defined by

\begin{eqnarray*}
&& L_m\cdot v=-\left(t^{m} D+\frac{m-2\epsilon}{2} t^{m}\right)v,
\\&&L_m\cdot\overline{v}=-\big(t^{m} D+m t^{m}\big)\cdot\overline{v}=-\overline{D t^{m}v},
\\&&G_p\cdot v=\big(t^{p-\epsilon} D\theta+
(p-\epsilon)t^{p-\epsilon}\theta\big)\cdot v=\overline{D t^{p-\epsilon} v},
\\&&G_p\cdot\overline{v}= -t^{p+\epsilon}v.
\end{eqnarray*}
We see that $(H_{\epsilon,\frac{1}{2}}(M))_{\bar0}\oplus D(H_{\epsilon,\frac{1}{2}}(M))_{\bar1}$ is a $\mathfrak{g}[\epsilon]$-submodule, namely, $M=D(M)$. For $k\in\Z,p\in\mathbb{Z}+\epsilon$  and $v\in M$, we have $$(G_pG_{k-p}-G_{p-1}G_{k-p+1})\cdot v=t^{k}v\in U(\mathfrak{g}[\epsilon])v,\quad G_p\cdot v=\overline{Dt^{p-\epsilon}v}\in U(\mathfrak{g}[\epsilon])v.$$
For any $i\in\mathbb{N}$ and $v\in M$, we have $(-L_0+\epsilon)^i\cdot v= D^iv\in U(\mathfrak{g}[\epsilon])v$.  Evidently, $(H_{\epsilon,\frac{1}{2}}(M))_{\bar 0}\oplus D  (H_{\epsilon,\frac{1}{2}}(M))_{\bar 1}$ is
an irreducible $\mathfrak{g}[\epsilon]$-submodule of $H_{\epsilon,\frac{1}{2}}(M)$.

{\rm(ii)} Now we define the following linear map
\begin{eqnarray*}
\phi: (H_{\epsilon,\frac{1}{2}}(M))_{\bar0}\oplus D(H_{\epsilon,\frac{1}{2}}(M))_{\bar1} &\longrightarrow& \Pi(H_{\epsilon,0}(M))
\\ -t^{2\epsilon}v&\longmapsto&\overline{v}
\\  \overline{Dv}&\longmapsto& v,
\end{eqnarray*}
where $v\in M,\overline{v}\in\overline{M}$.
For $m\in\Z$, one can check that
 $\phi$ is a  $\mathfrak{g}[\epsilon]$-module isomorphism.

If $H_{\epsilon,0}(M)$ is  reducible, from Case $1$ in Lemma \ref{lemm200} (i),  one gets $${H'_{\epsilon,0}(M)}=\big(\C[t^{\pm1}]\oplus\overline{\C[{t}^{\pm1}]}\big)/ \C.$$ We see that the
following linear map $\psi$ is a  $\mathfrak{g}[\epsilon]$-module isomorphism, where
\begin{eqnarray*}
\psi:\quad (H_{\epsilon,\frac{1}{2}}(M))_{\bar0}\oplus D(H_{\epsilon,\frac{1}{2}}(M))_{\bar1}\quad&\longrightarrow& {H'_{\epsilon,0}(M)}
\\-t^{m+2\epsilon}&\longmapsto&\overline{t^m}
\\  \overline{D t^n}&\longmapsto& t^n,
\end{eqnarray*}
for $ m\in\Z,n\in\Z^*.$

\end{proof}



Based on Proposition \ref{I} and Lemmas \ref{th36544433} -\ref{lemm201}, we can derive the following theorem.

\begin{theo} \label{main2}
Suppose that $M$ is an irreducible $\mathcal{D}$-module.
Then $H_{\epsilon,b}(M)$ is irreducible if and only if one of the following holds:
 \begin{itemize}
\item[{\rm (i)}] $b\notin\{0,\frac{1}{2}\};$
\item[{\rm (ii)}] $b=0$ and  $M\ncong\C[t^{\pm1}];$

\item[{\rm (iii)}] $b=\frac{1}{2}$ and $M=D(M)$.
 \end{itemize}

\end{theo}

\subsection{Natural isomorphisms}

\begin{lemm}\label{lemm20011}
Let $b_1, b_2\in \C, b_1\notin\{0,\frac{1}{2}\}$ and $M_1,M_2$ be irreducible $\mathcal{D}$-modules.
 Then  $H_{\epsilon, b_1}(M_1)\cong H_{\epsilon,b_2}(M_2)$, as   $\mathfrak{g}[\epsilon]$-modules,  if and only if
 $b_1=b_2$ and   $M_1\cong M_2$, as $\mathcal{D}$-modules.
\end{lemm}
\begin{proof}

 The sufficiency is evident, so we only need to provide a proof for the necessity. Suppose that $\psi: H_{\epsilon, b_1}(M_1) \rightarrow H_{\epsilon,b_2}(M_2)$ is a $\mathfrak{g}[\epsilon]$-module isomorphism. For any $p\in \mathbb{Z}+\epsilon$, $d\in\mathbb{Z}^*$, and $0\neq v\in M_1$, we have $\psi(T_{p,d}\cdot v)=T_{p,d}\cdot \psi(v)$ by applying Lemma \ref{lemm2.3}. This, together with the same lemma, yields
\begin{eqnarray}\label{eq401}
b_1(1-2b_1)\psi(t^{p-\epsilon}\theta v)=b_2(1-2b_2)t^{p-\epsilon}\theta \psi(v).
\end{eqnarray}
Note that $b_1\neq 0$ and $b_1\neq\frac{1}{2}$. By Lemma \ref{lemm3.1}, we see that the action of $\partial_\theta$ on $v\in M_1$ is trivial. Then, using the fact that
$$G_{-\epsilon}\cdot\psi \left(\frac1{b_2(1-2b_2)}{\overline{t^{p-\epsilon} v}}\right)=
\psi \left(
\frac1{b_2(1-2b_2)}G_{-\epsilon}\cdot{\overline{t^{p-\epsilon}v}}\right),$$
we get
\begin{eqnarray}\label{eq411}
\frac1{b_2(1-2b_2)}\psi({t^{p-\epsilon} v})=\frac1{b_1(1-2b_1)}{t^{p-\epsilon} \psi(v)}.
\end{eqnarray}
Setting $p=\epsilon$ in \eqref{eq411}, we obtain
\begin{eqnarray}\label{eq422}
b_1(1-2b_1)=b_2(1-2b_2).
\end{eqnarray}
Substituting \eqref{eq422} into \eqref{eq401} and \eqref{eq411}, we immediately obtain
$$\psi(\overline{t^{p-\epsilon}  v})=\overline{t^{p-\epsilon} \psi(v) }\quad \mathrm{and} \quad \psi(t^{p-\epsilon}  v)=t^{p-\epsilon} \psi(v),$$
where $p\in\mathbb{Z}+\epsilon$. For any $i\in\mathbb{N}$ and $v\in V$, we have
$$\psi((-L_0+2\epsilon b)^i\cdot v)=(-L_0+2\epsilon b)^i\cdot\psi(v),$$ which implies $\psi( D^iv)= D^i\psi(v)$. Therefore, we conclude that $M_1\cong M_2$.

Now, for any $m\in\mathbb{Z}$ and $0\neq v\in M_1$, from $\psi(L_m\cdot v)=L_m\cdot\psi(v)$, it is easy to check that $b_1=b_2$.
\end{proof}

\begin{lemm}\label{lemm201111}
Suppose  that $M_1$ and $M_2$ are irreducible $\mathcal{D}$-modules.
 Then
  $H_{\epsilon,0}(M_1)\cong H_{\epsilon,0}(M_2)$, as   $\mathfrak{g}[\epsilon]$-modules, if and only if
    $M_1\cong M_2$, as $\mathcal{D}$-modules.
\end{lemm}
\begin{proof}
The sufficiency of the conditions is clear. Suppose that $\psi:H_{\epsilon,0}(M_1)\rightarrow  H_{\epsilon,0}(M_2)$ is an  isomorphism. For any $\overline{v}\in \overline{M_1}$, we note that $(-L_0-\epsilon)\cdot\overline{v}= \overline{Dv}$.
We consider the action of $D-k$ on $\overline{v}$.

If $(D-k)\cdot\overline{v}=0$ for some $k\in\Z$ and a nonzero $\overline{v}\in \overline{M_1}$, then $ \overline{D(t^{-k}v)}=0$ where $t^{-k}v\neq0$.
According to  Case $1$ in the proof of Lemma \ref{lemm200}, we know that the $\W$-module $M_1\oplus \overline{M_1} \cong \C[t^{\pm1}]\oplus\overline{\C[{t}^{\pm1}]}$.
Similarly, using $(D-k)\psi(\overline{v})=0$ in $\overline{M_2}$, we deduce that $M_2\oplus\overline{M_2}\cong\C[t^{\pm1}]\oplus \overline{\C[{t}^{\pm1}]}$. Thus $M_1\cong M_2$.

Consider that $D-k$ is injective on both $\overline{M_1}$ and $\overline{M_2}$ for all $k\in\Z$. By  Lemma \ref{lemm3.1}, we see that any $v\in M_1$ can be annihilated by $\partial_\theta$.
For any $v\in M_1, p\in\Z+\epsilon$, we have
\begin{eqnarray}\label{908e}
\nonumber \psi(G_p\cdot v)&=&\psi(\overline{t^{p-\epsilon}Dv})=\psi(\overline{(D-p+\epsilon)t^{p-\epsilon}v})
\\&=&\psi((-L_0-p)\cdot\overline{t^{p-\epsilon}}v)=(D-p+\epsilon)\cdot\psi(\overline{t^{p-\epsilon}v}).
\end{eqnarray}
Inserting \eqref{908e} into
 $\psi(G_p\cdot v)=G_p\cdot\psi(v)$,
 it is easy to check  $\psi(\overline{t^{p-\epsilon}v})=\overline{t^{p-\epsilon}\psi(v)}$ for $p\in\Z+\epsilon$.
 Then by $\psi(G_{-\epsilon}\cdot(\overline{t^{p-\epsilon}v}))=G_{-\epsilon}\cdot(\overline{t^{p-\epsilon}\psi(v)})$,
we see that
$\psi(t^{p-\epsilon}v)=t^{p-\epsilon}\psi(v)$ for all $v\in M_1$ and $p\in\Z+\epsilon$.   For $i\in\N,v\in M_1$,  from   $\psi(L_0^i\cdot v)=L_0^i\cdot\psi(v)$,  one has  $\psi( D^iv)= D^i\psi(v)$.  Thus $M_1\cong M_2$
in this case.
\end{proof}

The isomorphism theorem can be obtained from Lemmas \ref{lemm201},  \ref{lemm20011} and  \ref{lemm201111}.

\begin{theo}\label{th4333}
Suppose that $b_1,b_2\in\C$, $M_1,M_2$ are irreducible $\mathcal{D}$-modules.
Then $H_{\epsilon,b_1}(M_1)\cong H_{\epsilon,b_2}(M_2)$ as   $\mathfrak{g}[\epsilon]$-modules if and only if  one of the following holds:
\begin{itemize}
    \item[{\rm (i)}] $b_1=b_2,$ $M_1\cong M_2;$
    \item[{\rm (ii)}] $(b_1,b_2)=(\frac{1}{2},0),$ $M_1\cong M_2$, and $M_1= D  (M_1);$
    \item[{\rm (iii)}] $(b_1,b_2)=(0,\frac{1}{2}),$ $M_1\cong M_2$, and $M_2= D  (M_2)$.

\end{itemize}

\end{theo}

\begin{coro}
Let $b_1, b_2\in \C$. Then $H_{\epsilon, b_1}\cong H_{\epsilon,b_2}$ if and only if  $b_1=b_2$.

\end{coro}




\section{Applications and   examples}\label{ex666}

In this section, we utilize the functors $H_{\epsilon,b}$ to recover several  known irreducible $\mathfrak{g}[\epsilon]$-modules. Furthermore, we construct lots of new irreducible $\mathfrak{g}[\epsilon]$-modules.

\subsection{Intermediate series modules}
Let $\alpha\in\C[t^{\pm1}]$. Set  $\tau= D-\alpha$ in Theorem \ref{211} (1),
we get the irreducible $\mathcal{D}$-module
$M_{\alpha}=\mathcal{D} /\mathcal{D}\tau$
with a basis $\{t^k\mid k\in\Z\}$, with
\begin{eqnarray*}
D\cdot t^n=t^n(\alpha+n),\quad t^m\cdot t^n=t^{m+n}, \ \forall m,n\in\Z.
\end{eqnarray*}

For $\epsilon\in\{0,\frac{1}{2}\}, b\in\C$,    we obtain    $\mathfrak{g}[\epsilon]$-module
$M_{\epsilon,\alpha,b}=H_{\epsilon,b}(\C[t^{\pm1}])=\C[t^{\pm1}]\oplus  \overline{\C[t^{\pm1}}]$
 with the following  actions:
\begin{eqnarray*}
 && L_m\cdot t^n=-(\alpha+n+(m-2\epsilon)b)t^{m+n},
 \\&&  L_m\cdot\overline{t^n}=-\left(\alpha+n+m(b+\frac{1}{2})+\epsilon(1-2b)\right)\overline{t^{m+n}},
  \\&&
 G_p\cdot t^n=(\alpha+n+2(p-\epsilon)b)\overline{t^{p-\epsilon+n}},
 \\&&  G_p \cdot\overline{t^n}=-t^{p+\epsilon+n},
 \end{eqnarray*}
 for $m,n\in\Z,p\in\mathbb{Z}+\epsilon$.

Especially,  if $\alpha\in\C$, then
 $M_{\epsilon,\alpha,b}$  is isomorphic to the intermediate series weight modules  of $\mathfrak{g}[\epsilon]$ in \cite{S}.  Moreover, we have

 \begin{coro}[\cite{S}]
 $M_{\epsilon,\alpha,b}$ is irreducible except that $\alpha\in\Z+\epsilon,b=0$ or $\alpha\in\Z,b=\frac{1}{2}$.
 \end{coro}

 \begin{rema}
Suppose that $\alpha\in\C[t^{\pm1}]\setminus\C$. It is clear that $M_{\epsilon,\alpha,b}$  is neither a weight module  nor a restricted modules for $\mathfrak{g}[\epsilon]$.
In particular, $\M_{\epsilon,\alpha,b}$ is  irreducible  if and only if $b\neq\frac{1}{2}$.

 \end{rema}

\subsection{$U(\mathfrak{h})$-free modules }
For $\lambda\in\C^*$, we introduce the irreducible $\mathcal{D}$-module $\Omega(\lambda)$ defined in Theorem \ref{211}, which has a basis $\{D^n\mid n\in\N\}$, and the $\mathcal{D}$-actions  are defined as follows
\begin{eqnarray*}
&&t^m\cdot (-D)^n=\lambda^m(-D+m)^n, \quad  D\cdot D^m= D^{m+1}
\end{eqnarray*}
for $m\in\Z, n\in\N$.

For   $b\in\C$,
 we have the $\mathfrak{g}[\epsilon]$-module  $\Omega_{\epsilon}(\lambda,b)=H_{\epsilon,b}(\Omega(\lambda))
$ with the  actions:
\begin{eqnarray*}
 && L_m\cdot (-D)^n=\lambda^m\big( -D+m(1-b)+2b\epsilon \big)(-D+m)^{n},\\
 && L_m\cdot\overline{(-D)^n}=\lambda^m\left( -\overline{D}+m(\frac{1}{2}-b)+(2b-1)\epsilon\right){(-\overline{D}+m)^{n}},\\
 &&G_p\cdot (-D)^n=-\lambda^{p-\epsilon}\left(-\overline{D}+2p(\frac{1}{2}-b)+(2b-1)\epsilon \right){(-\overline{D}+p-\epsilon)^{n}},\\
 && G_p\cdot\overline{(-D)^n}=-\lambda^{p+\epsilon}(-D+p+\epsilon)^{n},
 \end{eqnarray*}for $m\in\Z,p\in\mathbb{Z}+\epsilon, n\in\N$.

It can be seen that the $U(\frak h)$-module, which was constructed and analyzed in \cite{YYX1} through direct calculation, is isomorphic to $\Omega_\epsilon(\lambda,b)$. It follows from Theorem  \ref{main2} that

\begin{coro}[\cite{YYX1}]
The module $\Omega_\epsilon(\lambda,b)$ is irreducible
 if and only if $b\neq\frac{1}{2}$.
\end{coro}

\subsection{Degree two modules}\label{63ty}

Let $f(t)\in \C[t^{\pm1}]$ be such that $ D^2-f(t)\in\C(t)[ D]$ is an irreducible element. Set $\tau= D^2-f(t)$ in Theorem \ref{211} (1). Then we   obtain the irreducible $\mathcal{D}$-module
$M=\mathcal{D}/\big(\mathcal{D}\cap(\C(t)[ D]\tau)\big)$
with a basis $\{t^k,   t^k D, \mid  k\in\Z\}$. The $\mathcal{D}$-actions on $M$ are presented as
 \begin{eqnarray*}
 &&t^m\cdot t^n=t^{m+n},\quad t^m\cdot(t^n D)=t^{m+n} D,\quad
   D\cdot t^n= t^{n}(D+n), \quad
 D\cdot(t^n D)=t^{n}(f(t)+n D)
\end{eqnarray*}
for $m, n\in\Z$.

For $b\in\C$, we have  $\mathfrak{g}[\epsilon]$-module $H_{\epsilon,b}(M)$ with the following  actions
\begin{eqnarray*}
  &&L_m\cdot t^n=-t^{m+n}( D+n+(m-2\epsilon)b),
  \\&&
   L_m\cdot (t^n D)=-t^{m+n}(f(t)+(n+(m-2\epsilon)b) D),\\
 &&L_m\cdot (\overline{t^n})=  -\overline{{t}^{m+n}}(\overline{D}+n+m(b+\frac{1}{2})+\epsilon(1-2b)),
 \\&&
 L_m\cdot(\overline{t^n D}) =- \overline{t^{m+n}}(\overline{f(t)}+(n+m(b+\frac{1}{2})+(1-2b)\epsilon) \overline{D}),\\
 &&G_p\cdot t^n=  \overline{t^{p-\epsilon+n}}( \overline{D}+n+2(p-\epsilon)b),
 \\&&
 G_p\cdot ( t^n D)=  \overline{t^{p-\epsilon+n}}(\overline{f(t)}+(n+2(p-\epsilon)b) \overline{D}),\\
 &&G_p\cdot\overline{t^n}=-t^{p+\epsilon+n},\\  &&G_p\cdot\overline{t^n D}=-t^{p+\epsilon+n} D,
 \end{eqnarray*}
 where $m,n\in\Z, p\in\mathbb{Z}+\epsilon$.  From Theorems \ref{main2},
 we  see  that   the new  degree two modules $H_{\epsilon,b}(M)$ is irreducible   if  and only if $b\neq\frac{1}{2}$.

\subsection{Fraction modules}\label{65ty}
In this subsection, we will define a class of  new irreducible modules called  fraction modules over  $\mathfrak{g}[\epsilon]$.

Let $n\in\N,\alpha=(\alpha_0,\alpha_1,\ldots,\alpha_n)\in\C^{n+1},(b_0,b_1,\ldots,b_n)\in\C^{n+1}$
with $b_0=0$ and $b_i\neq b_j$ for all $i\neq j$. Set
$\tau=\frac{d}{dt}-\sum_{i=0}^n\frac{\alpha_i}{t-b_i}$ in Theorem  \ref{211} (1).
Then we get  the irreducible $\mathcal{D}$-module
\begin{eqnarray*}
M=\mathcal{D}/\big(\mathcal{D}\cap(\C(t)[ D]\tau)\big)
 \subset\C[t,(t-b_i)^{-1}\mid i=0,1,\ldots,n].
\end{eqnarray*}
 The actions of $\mathcal{D}$ on $M$ are defined as follows
\begin{eqnarray*}
&&\frac{d}{dt}\cdot f(t)=\frac{d}{dt} (f(t))+f(t)\sum_{i=0}^n\frac{\alpha_i}{t-b_i},
\quad
t^m\cdot f(t)=t^m f(t),
\quad \partial_\theta\cdot f(t)=0,
\end{eqnarray*} where
  $f\in M, m\in\Z.$

 For any $b\in\C$,
  we get   the actions of   $\mathfrak{g}[\epsilon]$
on $H_{\epsilon,b}(M)$ as follows:
\begin{eqnarray*}
  &&L_m\cdot f(t)=-t^{m+1}\frac{d}{dt} (f(t))-t^{m+1}f(t)\sum_{i=0}^n\frac{\alpha_i}{t-b_i}-(m-2\epsilon)b t^mf(t),\\
 &&L_m\cdot\overline{f(t)}=-\overline{t^{m+1}} \frac{d}{dt} \overline{f(t)}-\overline{t^{m+1}f(t)}\sum_{i=0}^n\frac{\alpha_i}{\overline{t}-b_i}-(m(b+\frac{1}{2})+\epsilon(1-2b)) \overline{t^m f(t)},\\
 &&G_p\cdot f(t)=\overline{t^{p-\epsilon+1}} \frac{d}{dt} \overline{f(t)}+\overline{t^{p-\epsilon+1}f(t)}\sum_{i=0}^n\frac{\alpha_i}{\overline{t}-b_i}+2(p-\epsilon)b \overline{t^{p-\epsilon} f(t)},\\
 && G_p\cdot\overline{f(t)}=-t^{p+\epsilon}f(t).
 \end{eqnarray*}
 By Theorems \ref{main2},
 we  know  that   the fraction module   $H_{\epsilon,b}(M)$ is  irreducible if and only if  $b\neq\frac{1}{2}$.





\subsection{Special  degree $n$ modules}\label{64ty}
For any $n\in\N$, some degree $n$ irreducible elements in $\C(t)[ D_t]$ were determined in  Lemma 16 of \cite{LZ}.

For any $n\in\N$, choose $\tau=\left(\frac{d}{dt}\right)^n-t$ in Theorem   \ref{211} (1). Thus, it is easy to get  the irreducible $\mathcal{D}$-module
$M=\mathcal{D}/\big(\mathcal{D}\cap(\C(t)[ D_t]\tau)\big)$
with a basis
$$\left\{t^k,\left(\frac{d}{dt}\right)^m\mid k\in\Z, m=0,1, \ldots, n-1\right\}.$$ The actions of $\mathcal{D}$ are defined as

\begin{eqnarray*}
  &&t^k\cdot \left(t^i\left(\frac{d}{dt}\right)^m\right)=t^{k+i}\left(\frac{d}{dt}\right)^m\quad \mathrm{for} \ k,i\in\Z, 0\leq m\leq n-1,\\
  && \frac{d}{dt}\cdot \left(t^i\left(\frac{d}{dt}\right)^m\right)=it^{i-1}\left(\frac{d}{dt}\right)^m+t^i\left(\frac{d}{dt}\right)^{m+1}\quad \mathrm{for} \ i\in\Z, 0\leq m< n-1,\\
 && \frac{d}{dt}\cdot \left(t^i\left(\frac{d}{dt}\right)^{n-1}\right)=it^{i-1}\left(\frac{d}{dt}\right)^{n-1}+t^{i+1}\quad \mathrm{for} \ i\in\Z.
 \end{eqnarray*}

For   $b\in\C$, we obtain the $\mathfrak{g}[\epsilon]$-module $H_{\epsilon,b}(M)$
with
\begin{eqnarray*}
  &&L_k\cdot \left(t^i\left(\frac{d}{dt}\right)^m\right)=-(i+b (k-2\epsilon))t^{k+i}\left(\frac{d}{dt}\right)^m-t^{k+i+1}\left(\frac{d}{dt}\right)^{m+1},\\
 &&L_k\cdot \left(t^i\left(\frac{d}{dt}\right)^{n-1}\right)=-(i+b (k-2\epsilon))t^{k+i}\left(\frac{d}{dt}\right)^{n-1}-t^{k+i+2},\\
 &&L_k\cdot \left(\overline{t^i\left(\frac{d}{dt}\right)^m}\right)=-\left(i+b(k-2\epsilon)+\frac{k+2\epsilon}{2}\right) \overline{t^{k+i}\left(\frac{d}{dt}\right)^m}-\overline{t^{k+i+1}\left(\frac{d}{dt}\right)^{m+1}},\\
&&L_k\cdot (\overline{t^i\left(\frac{d}{dt}\right)^{n-1}})=-\left(i+b(k-2\epsilon)+\frac{k+2\epsilon}{2}\right) \overline{t^{k+i}\left(\frac{d}{dt}\right)^{n-1}}-  \overline{t^{k+i+2}},
 \\&&
 G_p\cdot \left(t^i\left(\frac{d}{dt}\right)^m\right)=(i+2 (p-\epsilon)b)\overline{t^{p+i-\epsilon}\left(\frac{d}{dt}\right)^m}+\overline{t^{p+i+1-\epsilon}\left(\frac{d}{dt}\right)^{m+1}},
 \\
 &&G_p\cdot \left(t^i\left(\frac{d}{dt}\right)^{n-1}\right)=(i+2 (p-\epsilon)b)\overline{t^{p+i-\epsilon}\left(\frac{d}{dt}\right)^{n-1}}+\overline{t^{p+i+2-\epsilon}},
 \\
 &&G_p\cdot \left(\overline{ t^i \left(\frac{d}{dt}\right)^m}\right)=-t^{p+i+\epsilon}\left(\frac{d}{dt}\right)^m,
 \\&& G_p\cdot (\overline{t^i\left(\frac{d}{dt}\right)^{n-1}})=-t^{p+i+\epsilon}\left(\frac{d}{dt}\right)^{n-1},
 \end{eqnarray*}
  where  $k,i\in\Z,p\in\mathbb{Z}+\epsilon,0\leq m< n-1$.   By Theorems  \ref{main2},
 we  know  that   the special degree $n$ modules   $H_{\epsilon,b}(M)$ is  irreducible if and only if  $b\neq\frac{1}{2}$.
\begin{remark}
What we want to emphasize is that the degree two module, fractional module, and special degree $n$ module discussed in sections 5.3, 5.4, and 5.5 are novel, and to the best of our knowledge, no similar constructions have been found in the existing literature.
\end{remark}
\section*{Acknowledgements}
This work was carried out during the first author's visit to University of California, Santa Cruz. He would like to
thank Professor Chongying Dong  for hospitality during his visit. The authors also thank the anonymous referees for helpful
suggestions.
This work was supported by the National Natural Science Foundation of China (Grant Nos. 12171129, 12071405, 11971315).

\end{document}